\documentclass[12pt]{article}
\usepackage{amsmath,amsfonts, amssymb,verbatim}

\newtheorem{pkt}{}[section]  
\newcommand{\bpk}{\begin{pkt}\rm }  
\newcommand{\epk}{\end{pkt}}  
\newtheorem{Lemm}{Lemma}[subsection]  
   
\newtheorem{Prop}{Proposition} 
   
\newtheorem{Cor}{Corollary}

\newcommand{\Q}{{\mathbb Q}}

\newcommand{\C}{{\mathbb C}}

\newcommand{\be}{\begin{equation}}  
 \newcommand{\ee}{\end{equation}}   

\newcommand{\trd}{{\rm tr.d.}}   
\newcommand{\ld}{{\rm ldim}}

\newcommand{\bl}{\begin{Lemm}}   
   
\newcommand{\el}{\end{Lemm}}   
   
\newcommand{\bt}{\begin{Theo}}   
   
\newcommand{\et}{\end{Theo}}   
   
\newcommand{\bp}{\begin{Prop}}   
   
\newcommand{\ep}{\end{Prop}}   
   
\newcommand{\inv}{^{-1}}   
   
\newcommand{\bc}{\begin{Cor}}   
   
\newcommand{\ec}{\end{Cor}}

\newcommand{\ra}{\rangle}  
\newcommand{\la}{\langle}  
   
\newcommand{\subs}{\subseteq}

\newcommand{\cl}{\mathrm{cl}}

\begin{document}
\author{K.Tent and B.Zilber}
\date{ October 2, 2012}
\title{A non-desarguesian projective plane}

\maketitle
\abstract{We construct a new non-desarguesian projective plane  from a complex analytic  structure.  At the same time the construction can be explained in terms of so called Hrushovski's construction. This supports the hypothesis  that in general structures produced by Hrushovski's construction have ``prototypes'' in complex geometry. }

\section{Introduction}
Hrushovski's construction of ``new'' strongly minimal  structures and more generally ``new'' stable  structures proved very effective in providing a number of examples to 
classification problems in stability theory. For example, J.Baldwin used this method to construct a non-desarguesian projective plane of Morley rank 2 (see e.g. \cite{TZ}). But there is still a classification
problem of similar type which resists all attempt of solution, the Algebraicity (or Cherlin-Zilber) Conjecture. At present there is a growing belief that there must exists a simple 
group of finite   Morley rank which is not isomorphic to a group of the form $G(\mathbb{F})$ for $G$ an algebraic group and $\mathbb{F}$ an algebraically closed field (a {\em bad} group).

The second author developed an alternative interpretation of the ``new'' stable  structures obtained by Hrushovski's construction, see e.g. \cite{Zparis}. In this interpretation  
the universe $M$ of the structure is represented by a complex manifold and   relation by some subsets of $M^n$ explained in terms of the analytic structure on $M.$ In this interpretation
Hrushovski's predimension inequality corresponds to a form of (generalised) Schanuel's conjecture. We argue that looking for stable structures of analytic origin is potentially a better 
way  of producing new stable structures.

Below we briefly explain a construction of a new non-desarguesian projective plane that originates in a complex analytic  structure. The new, in comparison with previous examples of e.g.
``green fields'' (see \cite{Zbicol}) is that we have to use a non-trivial collapse procedure.
\section{An $\omega$-stable analytic structure}
\bpk \label{Liou0} Consider structures $K_f=(K,+,\cdot,f),$ where $(K,+,\cdot)$ is a field  and $f:K\to K$ a unary function.

Let $L^{alg}$ be a relational language for structures of the form $K_f,$
 relations of which are  those of $L$ along with all the relations corresponding to Zariski closed 0-definable subsets of $K^n.$ We always assume that $K$ is a field of characteristic $0.$
  Let
$C(K_f)$ be the  class of all finite $L^{alg}$-structures that can be embedded in $K_f.$
  
 Note that in the language $L^{alg}$ we can say for an $n$-tuple $X$ and a variety $W$ over $\Q$ that $X\in W.$ So the expression $\trd(X)=m$ means that $m$ is
 the dimension of the smallest variety $W$ over $\Q$ such that $X\in W.$ 

\epk
\bpk \label{Liou1}
{\bf Theorem.} (A.Wilkie, P.Koiran, \cite{W},\cite{Ko})  The structure $\C_f=(\C,+\cdot,f),$ where $f$ is an entire {\bf Liouville function,} 
is a model of the first order $\omega$-stable theory $T_f.$  
For every finite subset $X$ of the structure holds the Hrushovski inequality 
$$\delta(X)\ge 0,\mbox{ where } \delta(X):=\trd(X\cup f(X))-|X|.$$
\epk

\bpk \label{Liou2} {\bf Theorem.} $\C_f$ is $\omega$-saturated.

{\bf Proof.} Let $d$ be the dimension in $\C_f$ corresponding to the predimension $\delta.$
For $X\subs \C$ we denote $\cl(X)$ the closure of $X,$ that is the set of those $y\in \C$ that
satisfy $d(y/X)=0.$  We will say that $\C_f$ has CCP, the countable closure property, if
$\cl(X)$ is countable for any countable $X.$  

We need the following.
\epk
\bpk \label{Ccp}{\bf Lemma.} 
%(i) $\C_f$ is prime over a generic indiscernible subset  $I\subset \C;$
 $\C_f$ is the unique model of $T_f$ of cardinality continuum which satisfies the countable closure property.

{\bf Proof.} We use the general theory of ``pseudo-analytic structures'' developed in \cite{Zexp} in application to pseudo-exponentiation and in many further variations by other authors.  

Claim 1. $\C_f$ satisfies CCP.

This can be proved following lines of the proof of the same statement for $\C_{\exp}$ (Lemma 5.12) in \cite{Zexp},  only simpler. We do not  need to use Ax's theorem on Schanuel's property for function fields since the ``Schanuel'' property for $\C_f$ holds (by Wilkie's theorem). Also, instead of linear dependence one uses the equality relation, since in our $\delta$ we count the size of a set rather than linear dimension.

Now we recall the {\bf strong existential closedness axiom} of \cite{Zexp} (s.e.s.c) in the setting  of our structure.  It requires that for every algebraic variety $V\subs \C^{2n}$ over a subfield $k$ which satisfies the property that it is normal (rotund, in more recent terminology) and free with respect to =, there is a generic over $k$ point in $V$ of the form
$\la x_1,\ldots,x_n,f(x_1),\ldots,f(x_n)\ra.$ 

Claim 2. Any model of $T_f$ satisfies the   strong existential closedness axiom with respect to $f.$

Proof. If we omit the requirement of the point being generic, then we have   the   existential closedness axiom,  one of the axioms of $T_f.$ We claim that this is enough to have a generic point. Indeed, we can assume that $V$ is defined over
some finite self-sufficient $A\le \C$ and that $\dim V$ is minimal, i.e. equal to $n.$ Then
either $\delta(x_1,\ldots,x_n/A)=0$ for the point $\la x_1,\ldots,x_n,f(x_1),\ldots,f(x_n)\ra\in V,$ 
or $x_i=x_j$ for some $i\neq j,$ or $x_i=a\in  A$ for some $i.$ By existential closedness we still can find a point for which no of the equalities holds, so 
$$\trd(x_1,\ldots,x_n,f(x_1),\ldots,f(x_n)/A)=n=\dim V.$$
This proves that the point is generic. Claim proved.

The main theorem of \cite{Zexp} states that there is unique structure of a given uncountable cardinality that satisfies the Schanuel property,  strong existential closedness and CCP. Adapted to our setting, we derive using the claims that $\C_f$ is the unique structure of cardinality continuum satisfying the properties. $\Box$ Lemma.

Moreover, by  \cite{Zexp} and \cite{Zaec}
the unique model is $\omega$-homogeneous over submodels, so is $\omega$-homogeneous. Finally, by construction of the unique model, it is universal for the class of finitely generated substructures of the form $K_f.$ It follows that $\C_f$ is $\omega$-saturated. This complets the proof of the theorem. $\Box$

\epk
\bpk On class $C(\C_f)$ we define an equivalent predimension $\delta_0$ as follows.
For $Y\in C(\C_f)$
$$\delta_0(Y):= \trd(Y)-|Y^2\cap F|,$$
where $F$ is the graph of $f.$

\medskip

{\bf Lemma.} Let $M$ be a finite $L^{alg}$-structure. Then 
$\delta_0(Y)\ge 0$ for all $Y\subs M$ if and only if $M\in C(\C_f).$

{\bf Proof.} Suppose $M\in C(\C_f)$ and $Y\subs M.$ 
 Let $$Y^2\cap F=\{ \la x_i,f(x_i)\ra: i=1,\ldots n\}.$$
We assume the $x_i$ all distinct. Set $X=\{ x_1,\ldots,x_n\}.$ Then
by assumption $$0\le \delta(X)=\trd(X)-n\le \trd(Y)-|Y^2\cap F|=\delta_0(Y).$$
Conversely, suppose $\delta_0(Y)\ge 0$ for all $Y\subs M.$ Then for every
$X=\{ x_1,\ldots,x_n\}\subs M$ such that $\{ f(x_1),\ldots,f(x_n)\}\subs M,$ we have 
$\delta(X)\ge \delta_0(X\cup f(X))\ge 0.$ Extend $M$ to $M'$ by adding new elements
of the form $f(x)$ for every $x\in M$ such that $f(x)\notin M.$ We define such $f(x)$
to be mutually algebraically independent over $M.$  It is easy to see that for every
$X\subs M$ in regards to $f$ on $M'$ the inequality $\delta(X)\ge 0$ holds.
It follows that the diagram of $M'$ is consistent with the theory of $\C_f.$ By
\ref{Liou2} $M'$ can be embedded in $\C_f.$ $\Box$ 

\epk
\bpk \label{zero} {\bf Assumption $f(0)=0.$}

This can be achieved by setting $f(x)$ to be $f(x)-f(0).$ This does not effect the statement \ref{Liou1}, except for the change in the definition of the predimension $\delta,$ we have to replace it by the
predimension over $0.$ This does not effect our calculations with predimension below.

In particular,
if $x_1,\ldots,x_n$ is a generic (in the sense of Morley rank) tuple in the field $\C_f,$ then
\be
\label{eq} \delta(x_1,\ldots,x_n)=n\ee 
\epk
\section{Mild Collapse}

\bpk \label{alpha} Consider the class $C(\C_f).$ This is an amalgamation class with respect to strong embeddings $\le$ determined by $\delta.$

Let $\mu$ be a Hrushovski function satisfying $\mu(\alpha)=1,$ for any $\alpha,$ 
which is a code of a pair $(x,x_1,y_1,x_2,y_2/a_1,a_2,b)$ 
in a substructure 
$\{ x,x_1,y_1,x_2,y_2,a_1,a_2,b\}$
that satisfies relations
$$\begin{array}{lll}a_1x=x_1,\ a_2x=x_2,\\
f(x_1)=y_1,\ f(x_2)=y_2,\\
y_1-y_2=b
\end{array}$$

Note, that the code of type $\alpha$ says that $f(a_1x)-f(a_2x)=b,$ and $\mu(\alpha)=1$ amounts to saying that the latter has at most one solution in $x.$\\
%\epk
%\bpk {\bf Remark.} The condition \ref{alpha} on $\mu$ is consistent with setting $\mu(\beta)$ finite.

Consider the corresponding subclass $C_\mu(\C_f).$ We want to porve that this class has AP with respect to $\le.$
%and let $K_f$ be a countable  rich structure for this class.
\epk
\bpk \label{main}{\bf Proposition.}
 Let $B,M,N\in C_\mu(\C_f),$ $B\le M,\ B\le N,$ and let $M\otimes_B N\subs \C_f$ be a free amalgam over $B$ in $C(\C_f).$

 Suppose a code of type $\alpha$ is realised in $M\otimes_B N$ by $\{ x,x_1,y_1,x_2,y_2,a_1,a_2,b\}$ and $b\in M.$
Then 

(i) either $\{ x,x_1,y_1,x_2,y_2,a_1,a_2,b\}\subs M,$ 

(ii) or $b\in B$ and $\{ x,x_1,y_1,x_2,y_2,a_1,a_2,b\}\subs N,$

(iii) or $\{ x_1,x_2\}\subs B,$
$\{ y_1,y_2,b\}\subs M$ and $\{ x,a_1,a_2\}\subs N-M,$

(iv) or $\{ x_1,x_2,b\}\subs B,$
$\{ y_1,y_2,\}\subs N$ and $\{ x,a_1,a_2\}\subs M-N.$

{\bf Proof.} 
 Let $\Q(M),\Q(N)$ and $\Q(B)$ be the field generated by the corresponding subsets. Note, that since $B\le M,$ there is no new relations $f(u)=v$ on $\Q(M)$ and similarly with $\Q(N).$ 

Claim. $\Q(M)$ and $\Q(N)$  are linearly dijoint over $\Q({B}).$  
$M\otimes_B N$ can be naturally embedded in the free composite of fields $\Q(M)\otimes_{\Q(B)} \Q( N).$ 

The first follows from the freenes. The rest by strong embeddings. $\Box$ Claim.

We may assume that $1\in B$ and $B=\Q(B)\cap M= \Q(B)\cap N.$

We continue with auxiliary lemmas.

\epk
\bpk \label{yy} {\bf Lemma.} $\{ y_1,y_2,b\}\subs M,$ or $b\in B$ and  $\{ y_1,y_2,b\}\subs N.$ 

Claim. $\{ y_1,y_2,b\}\subs \Q(M),$ or $b\in B$ and  $\{ y_1,y_2,b\}\subs \Q(N).$ 

Case 1. $y_1,y_2\in N-M.$  By disjointness $\ld_{\Q(B)}(y_1,y_2,1)=\ld_{\Q(Bb)}(y_1,y_2,1).$ Since $b=y_1-y_2,$ it follows 
$b\in \Q(B),$ so $b\in B.$  

Case 2. $y_1\in N-M,$ $y_2\in M.$     Since $b=y_1-y_2,$ we have  $y_1\in \Q(M).$
$\Box$ Claim.

Now note that $\Q(M)\cap N=B,$ and Lemma follows by the Claim. $\Box$

\medskip

Without loss of generality we will assume below that  $\{ y_1,y_2,b\}\subs M.$

\epk
\bpk \label{xx}{\bf Lemma.} Assume $\{ y_1,y_2,b\}\subs M.$ Then
$\{ x_1,x_2\}\subs M.$

{\bf Proof.}  By freeness there are no relations of the form $f(u)=v$ between $N-B$ and $M-B.$ $\Box$
\epk
\bpk \label{aa} {\bf Lemma.} Assume $\{ y_1,y_2,b\}\subs M.$ Then
$\{ x,a_1,a_2\}\subs M,$ or $\{ x,a_1,a_2\}\subs N-M$ and $\{ x_1,x_2\}\subs B.$

{\bf Proof.} By definition in \ref{alpha} $x_1a_2=x_2a_1.$ 

Now, if $x\in M,$ then $a_1\in M,$ since $a_1=\frac{x_1}{x}.$ For the same reason $a_2\in M.$ So $\{ x,a_1,a_2\}\subs M.$

Similarly, if $a_1\in M,$  then $x=\frac{x_1}{a_1}\in \Q(M),$ $x\in M,$ and $\{ x,a_1,a_2\}\subs M.$

Hence, the alternative to  $\{ x,a_1,a_2\}\subs M$ is $\{ x,a_1,a_2\}\subs N,$
which implies $x_1,x_2\in N,$ since $x_1=a_1x$ and $x_2=a_2x.$ $\Box$
\epk
\bpk {\bf Lemma.} Assume $\{ y_1,y_2,b\}$ not a subset of $M.$ Then$\{ x,x_1,y_1,x_2,y_2,a_1,a_2,b\}$
$\{ y_1,y_2,b\}\subs N,$   and either $\{ x,a_1,a_2\}\subs N,$ or $\{ x,a_1,a_2\}\subs M-N$ and $\{ x_1,x_2\}\subs B.$

{\bf Proof.} This is just the symmetric case with proofs corresponding to that of Lemmas \ref{xx} and \ref{aa}. $\Box$

\medskip
 
This lemma completes the proof of Proposition~\ref{main}. $\Box$ 
\epk
\bpk \label{mainL} {\bf Lemma.} In \ref{main},
suppose $N$ or $M$ is minimal over $B.$ Then (i) or (ii) of \ref{main}  holds.

{\bf Proof.} Suppose $N$ is minimal.  Then 
(iii) is not possible, since each of $a_1,a_2$ and $x$ is algebraic over $M.$   

Under the same assumption, if (iv) holds, then  $y_1\in B$ or $y_2\in B,$ because $\delta_0(y_i/B)=0$ for $i=1,2.$ But then
both will have to be in $B$ since $y_1-y_2=b$ and $b\in B.$ This brings us into case (i).

Now consider the case   $M$ is minimal. Suppose (iii) holds. Then at least two of the elements of $\{ y_1,y_2,b\}$ has to be in $B,$ and again,
since $y_1-y_2=b,$ all three are in $B.$ This brings us into the case (ii).

(iv) can not hold since $\{ x,a_1,a_2\}\subs M-N$ is in contradiction with minimality of $M.$ $\Box$ 
\epk
\bpk {\bf Proposition.} $C_\mu(\C_f)$ is an amalgamation class.

{\bf Proof.} We consider $B\le M,$ $B\le N$ with an assumption that $b\in M$ and one of the extensions is minimal. Suppose  $M\otimes_B N$ is not in $C_\mu(\C_f),$ that is there are
$\{ x,x_1,y_1,x_2,y_2,a_1,a_2,b\}$ and $\{ x',x'_1,y'_1,x'_2,y'_2,a_1,a_2,b\},$  substructures of code $\alpha$
in  $M\otimes_B N$ such that $x\neq x'.$ 

By \ref{mainL} we will have $\{ x,x_1,y_1,x_2,y_2,a_1,a_2,b\}\subs M$ and $\{ x',x'_1,y'_1,x'_2,y'_2,a_1,a_2,b\}\subs N.$
Hence $\{ a_1,a_2,b\}\subs M\cap N=B.$ Now we may assume that $N$ is minimal over $B,$ that is $N=B\cup \{ x',x'_1,y'_1,x'_2,y'_2\}.$
Since the type of  $\{ x',x'_1,y'_1,x'_2,y'_2\}$ over $B$ given by code $\alpha$ is complete, we can identify  $\{ x',x'_1,y'_1,x'_2,y'_2\}$ with $\{ x,x_1,y_1,x_2,y_2\}$ thus identifying $M$ as an amalgam of $M$ and $N$ over $B.$  $\Box$
\epk
\bpk \label{Th1} {\bf Theorem.} There exists a countable rich structure $K_f$ for class $C_\mu(\C_f).$

(i) $K_f$ is an algebraically closed field with a function $f.$

(ii) Given $a_1\neq a_2$ and $b$ in $K_f,$ there is a unique solution to the equation $$f(a_1x)-f(a_2x)=b.$$
In particular, $f$ is a bijection on $K.$

  (iii) $K_f$ is embeddable in $\C_f.$  
  
(iv) depending on $\mu,$ the theory of $K_f$ is  $\omega$-stable of rank $\omega$ or strongly minimal. 

{\bf Proof.} (i) and (iv) follows from general theory.

(iii) is by \ref{Liou2}.

(ii) follows from the definition of $\mu.$ $\Box$ 
\epk
\section{Plane}
\bpk
Recall \cite{AS} that a {\bf ternary ring} $R$ is a set $R$ with two distinguished elements $0, 1$ and a ternary operation         
$T : R^3  \to R$ satisfying the following conditions:         
   
   (T1) $T (1, a, 0) = T (a, 1, 0) = a$ for all $a \in R;$
   
   (T2) $T (a, 0, c) = T (0, a, c) = c$ for all $a, c \in R;$
     
   (T3) If $a, b, c \in R,$ the equation $T (a, b, y) = c$     
         has a unique solution $y;$
                                 
   (T4) If $a, a' , b, b' \in R$ and $a \neq a',$ the equations  
         $T (x, a, b) = T (x, a' , b' )$ have a unique solution $x$ in $R;$
                                           
   (T5) If $a, a' , b, b' \in R$ and $a \neq a',$ the equations  
   
         $T (a, x, y) = b,\ T (a' , x, y) = b'$ have a      
         unique solution $x, y$ in $R.$  \\
            
Consider the ternary operation $$T(a,x,b)=f\inv(f(ax)+b)$$ on $K_f.$
\epk
\bpk \label{f0} 
{\bf Lemma.} The ternary operation $T(a,x,b)$ on $K_f$ determine a ternary ring with $0$ and $1$ of the field $K.$.

{\bf Proof.} Check using \ref{Th1} and \ref{zero}. $\Box$
\epk
\bpk \label{AS} {\bf Theorem} (see \cite{AS}) Every projective plane $\mathbf{P}$ is bi-interpretable with a ternary ring $R$ (associated ternary ring). 

 Every Desarguesian plane has a unique
associated ternary ring, which is an associative division ring.
\epk
\bpk {\bf Corollary.} The projective plane associated with the ternary ring $(K_f,T)$ is not desarguesian. 

Proof. Suppose the projective plane is Desarguesian. Then by \ref{AS} the ternary ring $(K_f,T)$ is an associative division ring. 
In particular, we will have the identity $$T(a,x,b)=a*x \dot{+} b$$
for $a*x:= T(a,x,0)=ax$ and $x\dot{+} b:=T(1,x,b)=f\inv(f(x)+b).$

The associativity law will give us the identity
$$a(x_1+x_2)=ax_1+ax_2,\mbox{ equivalently } af\inv(f(x_1)+x_2)=f\inv(f(ax_1)+ax_2).$$
The latter identity implies $\delta(a,x_1,x_2)<3$ for any elements $a,x_1,x_2\in K_f,$ in  
contradiction with (\ref{eq}). $\Box$

\epk
                                            
\end{document}